\input amstex
\documentstyle{amsppt}
\document
\magnification=1200
\NoBlackBoxes
\nologo
\vsize18cm


\bigskip

\centerline{\bf THETA FUNCTIONS, QUANTUM TORI}

\medskip

\centerline{\bf AND HEISENBERG GROUPS}

\medskip

\centerline{\bf Yuri I. Manin}

\medskip

\centerline{\it Max--Planck--Institut f\"ur Mathematik, Bonn}

\bigskip

{\bf Abstract.} A linear algebraic group $G$
is represented by the linear space of its algebraic
functions $F(G)$ endowed with multiplication and comultiplication
which turn it into a Hopf algebra. Supplying $G$ with a Poisson
structure, we get a quantized version $F_q(G)$
which has the same linear structure and comultiplication,
but deformed multiplication. This paper develops 
a similar theory for abelian varieties.
A description of abelian varieties $A$ in terms of linear algebra
data was given by Mumford: $F(G)$ is replaced by the graded ring of
theta functions with symmetric automorphy factors,
and comultiplication is replaced by the
Mumford morphism $M^*$ acting on pairs of points as $M(x,y)=M(x+y,x-y).$
After supplementing this by a Poisson structure and replacing the classical theta functions
by the quantized ones, introduced by the author earlier, we obtain a structure
which essentially coincides with the classical one
so far as comultiplication is concerned, but has a deformed 
multiplication which moreover becomes only partial.  The classical graded ring 
is thus replaced by a linear category. Another important difference
from the linear case is that abelian varieties with different
period groups (for multiplication) and different quantization
parameters (for comultiplication) become interconnected
after quantization.

\bigskip

\centerline{\bf 0. Introduction and summary}

\medskip

{\bf 0.1. Classical theta functions.} Let $\bold{C}^d$ be a complex
vector space, $\Lambda^{\prime}\oplus \Lambda \subset \bold{C}^d$ a discrete sublattice of rank $2d$ split into the sum of two sublattices of rank $d$.
Classical theta functions are entire functions $\Theta (z), z\in \bold{C}^d$,
such that
$$
\Theta (z+\lambda^{\prime})=\Theta (z)\quad \roman{for\ all}\quad 
\lambda^{\prime} \in \Lambda^{\prime},
\eqno(0.1)
$$
$$
\Theta (z+\lambda ) = c(\lambda )\,e^{q(\lambda ,z)}\,\Theta (z)
\quad \roman{for\ all}\quad 
\lambda \in \Lambda,
\eqno(0.2)
$$
where $c:\,\Lambda \to\bold{C}$ is a map and $q:\, \Lambda\times \bold{C}\to
\bold{C}$ is a biadditive pairing linear in $z$.

\smallskip

$\Lambda^{\prime}\oplus\Lambda$ is called {\it the period lattice} for $\Theta$, and $(c,q)$ {\it the automorphy factors.} These data are not arbitrary.
It is convenient to break the classical restrictions, imposed on them,
into two groups.

\smallskip

(i) {\it Group action and Riemann symmetry conditions.} For each $\lambda^{\prime}+\lambda 
\in \Lambda^{\prime}\oplus\Lambda$, consider
the following map of the (total space of the) trivial line bundle 
$\bold{C}^d\times \bold{C}$ extending the shift by $\lambda^{\prime}+\lambda$
on the base space:
$$
t_{\lambda^{\prime}+\lambda}:\, (z, u)\mapsto (z+\lambda^{\prime}+\lambda , c(\lambda )\,e^{q(\lambda ,z)}\,u).
\eqno(0.3)
$$
Then we must have 
$$
t_{\lambda^{\prime}_1+\lambda^{\prime}_2+\lambda_1+\lambda_2}=
t_{\lambda^{\prime}_1+\lambda_1}\circ t_{\lambda^{\prime}_2+\lambda_2}.
\eqno(0.4)
$$
Solvability of (0.4) is essentially ensured by {\it Riemann's symmetry conditions}
imposed upon the lattice and the automorphy factors.
\smallskip

If they  hold, the quotient of $\bold{C}^d\times \bold{C}$
modulo the action of $\Lambda^{\prime}\oplus\Lambda$ thus defined 
will be a line bundle $L$ over the complex torus $\Cal{A}:=\bold{C}^d/(\Lambda^{\prime}\oplus\Lambda )$, and any theta
function satisfying (0.1) and (0.2) will descend to
a holomorphic section of $L$.

\smallskip

(ii) {\it Formal theta functions and Riemann's positivity conditions.}
Equations (0.1) are satisfied for the basic exponents
$e^{2\pi i\mu (z)}$ where $\mu$ runs over all linear functions
$\bold{C}^d\to \bold{C}$ taking integral values on
$\Lambda^{\prime}$. Such $\mu$ form an abstract lattice
which we will denote $H$. Usually one passes to a basis
of $\bold{C}^d$ which is simultaneously a basis
for $\Lambda^{\prime}$; then the basic exponents
are simply $e^{2\pi i\bold{n}^t\cdot z}$, $\bold{n}\in\bold{Z}^d.$ 

\smallskip

Any entire function satisfying
(0.1) can be expanded into a Fourier series
$$
\Theta (z) = \sum_{\mu\in H} a_{\mu}\,e^{2\pi i\mu (z)},
\eqno(0.5)
$$
whose coefficients decay swiftly enough. 

\smallskip

Equations (0.2) written for (0.5) 
translate into some linear recurrent relations
for coefficients $a_{\mu}$ which may have solutions satisfying
only weaker decay conditions, or even defining
everywhere divergent series. We will interpret all these
solutions as defining generally a linear space
of {\it formal} theta functions.

\smallskip

We can ask, when all formal solutions are in fact entire.
The answer is given by {\it Riemann's positivity conditions}.
(We will reproduce both conditions explicitly
below in the form more appropriate for our purposes).

\smallskip

If both symmetry and positivity conditions are satisfied
for the period lattice $\Lambda^{\prime}\oplus\Lambda$  
and the system of automorphy factors $(c,q)$,
the complex torus $\Cal{A}:=\bold{C}^d/(\Lambda^{\prime}\oplus\Lambda )$
has the canonical structure of {\it an abelian variety},
the line bundle $L$ on $\Cal{A}$ defined above is ample,
and the space of all theta functions satisfying (0.1), (0.2)
can be canonically identified with $\Gamma (\Cal{A},L).$

\smallskip

However, we will be definitely interested in more general
thetas.

\medskip

{\bf 0.2. Theta functins on an algebraic torus.} The total group
of shifts by $\Lambda^{\prime}\oplus\Lambda$
acting upon $\bold{C}^d$ directly and upon $\bold{C}^d\times\bold{C}$
via automorphy factors has the subgroup $\Lambda^{\prime}$.
Let us make an inventory of the data which we obtain after moding out
only this subgroup.

\smallskip

(a) {\it Algebraic torus $T(H)$ and its characters.}
This is an affine commutative algebraic group over $\bold{C}$ isomorphic
to a product of $d$ multiplicative groups $\bold{G}_m$. More precisely,
it is the spectrum of the group ring of
$H$, the abstract
lattice $\roman{Hom}\,(\Lambda^{\prime},\bold{Z})$ introduced
in (0.1) (ii). In the group ring, we write
the elements of $H$ multiplicatively, and then denote
them $e(h)$: 
$$
 e(g+h)=e(g)\,e(h),\quad g,h\in H.
\eqno(0.6)
$$ 
We have an analytic map, which is in fact an isomorphism:
$$
p:\,\bold{C}^d\to T(H)(\bold{C}):\quad p^*(e(h)):=e^{2\pi i\mu (.)},
\eqno(0.7)
$$
where $\mu (.)$ is the linear function on $\bold{C}^d$
extending $h$ as a function on $\Lambda^{\prime}.$

\smallskip

(b) {\it Multiplicative period lattice $B\subset T(H)(\bold{C})$.} By definition, this is the image of $\Lambda$ with respect to the map (0.7),
$B=p(\Lambda )$, so that the complex torus $\Cal{A}$ 
is now represented as
$$
\Cal{A} = T(H)/B 
\eqno(0.8)
$$

\smallskip

(c) {\it Automorphy factors and theta functions on $T$.}
Any formal theta function (0.5) can be interpreted as a formal
function on $T$, an infinite linear combination $\theta =\sum_{h\in H}a_h\,e(h).$
This takes care of equations (0.1).

\smallskip

Clearly, $p$ induces an isomorphism of $\Lambda$ with $B$.

\smallskip

After replacing $\lambda\in \Lambda$ by $b:=p(\lambda )\in B$,
the left hand side of (0.2) becomes the multiplicatively shifted series
which we now denote
$$
b^*(\theta ):= \sum_{h\in H} a_h\,h(b)\,e(h)
\eqno(0.9)
$$
where  $h(b)$ is the value of the character $e(h)$
at the point $b\in T(H)(\bold{C}).$ 

\smallskip

To treat (0.2), we simply replace $c(\lambda )$ by $c_b^{-1},\,b\in B$,
and $e^{q(\lambda ,z)}$ by $e(-h_b)$
in the definition of automorphy factors. Here $b\mapsto h_b$ 
is a certain map $B\to H$ (which is easily seen to be a group homomorphism,
if we want (0.4) to hold).  Extra inversions were introduced
in order to rewrite 
the equations (0.2) in the following form:
$$
\forall\, b\in B,\quad c_b\,e(h_b)\,b^*(\theta ) =\theta .
\eqno(0.10) 
$$
Of course, in the case $d=1$, $B=\{q^{\bold{Z}}\}\subset \bold{C}^*$,
we get Jacobi's elliptic theta functions, for example,
the basic theta
$$
\theta_q(t) =\sum_{n\in Z} q^{n^2}\,t^n
\eqno(0.11)
$$
where we write $t=e(h_0)$, $h_0$ being a generator of $H$.
It converges everywhere on $\bold{C}^*$ for $|q|<1$,
defines a distribution on $|t|=1$ for $|q|=1$
and is a formal theta--function generally, satisfying the equations
$$
q^{m^2}\,t^{m}\,\theta_q (q^{2m}t)= \theta_q (t)
\eqno(0.12)
$$
for all $m\in \bold{Z}.$

\smallskip

This reformulation of basic definitions gives some immediate benefits,
the most important of which is the possibility
to extend the formal theory to arbitrary base fields $K$,
and analytic theory to $p$--adic fields and, more generally,
to arbitrary complete discretely normed fields  instead of $\bold{C}$.
In fact, over such fields reasonable functions with
additive period lattices do not exist so that the equations (0.1)
have no interesting solutions. Moding them out,
we are left with (0.10), which do admit
analytic solutions and, in particular,
provide the theory of $p$--adic uniformization of abelian varieties
with multiplicative reduction (Tate, Morikawa).

\medskip

{\bf 0.3. Homogeneous coordinate ring of an abelian
variety.} Briefly, this theory proceeds as follows.
Consider an algebraic torus $T(H)$ over 
a complete normed field $K$, endowed with a period subgroup 
$B\subset T(H,1)(K)$
which is free of the same rank as $H$ and discrete.
Consider moreover all theta functions on $T(K)$
with multiplicative periods $B$ and automorphy factors
satisfying Riemann's symmetry and positivity conditions
(the latter can be defined also in the $p$--adic case, 
cf. the main text below). Assume that this set
(consisting of ``ample'' automorphic factors)
is non--empty. Denote by $\Gamma (\Cal{L})$ the space
of theta functions with automorphy factors $\Cal{L}.$

\smallskip

The set of automorphy factors is endowed with
the natural commutative and associative multiplication.
This holds even if we drop the positivity conditions.
However, with positivity conditions, all relevant
theta functions have swiftly decaying coefficients
and therefore as well can be multiplied in the space of
formal series in exponents. Moreover, a product of
theta functions with the same periods $B$ is again
a theta function with these periods.

\smallskip

This multiplication 
induces a mapping $\Gamma (\Cal{L}^{\prime})\otimes 
\Gamma (\Cal{L})\to \Gamma (\Cal{L}^{\prime}\Cal{L}).$

\smallskip

Construct the ring $\oplus \Gamma (\Cal{L})$ graded by
the semigroup of ample and symmetric $\Cal{L}$
(symmetric means invariant with respect to the involution
$x\mapsto x^{-1}$ of $T(H,1)$).
Its projective spectrum is an abelian variety $A$
which is a model of 
$T(H,1)/B$ in algebraic geometry. Automorphy factors $\Cal{L}$ can be interpreted then
as invertible sheaves on $A.$

\smallskip

This beautiful theory still lacks some
important structures present over $\bold{C}$. Name\-ly, consider 
an abelian variety $\Cal{A}$ over the complex field.
Two theta functions on the {\it universal cover}
$\varphi :\,\widetilde{\Cal{A}}\to\Cal{A}(\bold{C})$
(as opposed to the toric cover $T(H)(\bold{C})$), 
 may look differently
but correspond to the same section of the same line
bundle $L$ on $\Cal{A}$, because we can choose
different splittings of the period lattice
$\roman{Ker}\,\varphi$ and change the trivialization
of $\varphi^*(L)$ appropriately. 

\smallskip

In the analytical language, this translates into the
functional equations for theta functions expressing
their modular properties with respect to the matrix fractional linear transformation acting upon both
$z$ and the half of the period matrix corresponding to $\Lambda$.

\smallskip

We cannot expect  a simple--minded toric version
of these equations to exist in general. For example, in the
$p$--adic case an additive uniformization of abelian
varieties does not exist because free abeian subgroups
of $p$--adic vector spaces cannot be discrete.

\smallskip

Notice that this sort of difficulties in principle 
vanishes in non--commutative geometry, where ``bad quotient spaces'' like
$\bold{Q}_p^d/\bold{Z}^d$ might have a perfectly well defined and
rich theory. Moreover, in non--commutative geometry we can
form another model of $T(H)/B$, independently of
whether $B$ is discrete or not. Such a model is represented
by the cross product of a ring of functions of $T(H)$
with the transformation group $B$, which is a special case
of {\it quantum torus}.

\smallskip

Remarkably, over $\bold{C}$ a replacement of modularity
for some quantum tori with maximally non--discrete (unitary) 
periods does exist and is 
expressed in the form of a beautiful
Morita equivalence: see [RiSch] and the references therein.

\smallskip

We return now to the algebraic tori. 

\medskip

{\bf 0.4. Heisenberg groups.} In the toric context, over an arbitrary
base field $K$, we can also
succintly express the conditions (0.4). To this end
let us introduce {\it the Heisenberg group} $\Cal{G}(H)$
of the torus $T(H)$. By definition, this is the group
of {\it linear endomorphisms} of the space of
(formal or algebraic) functions on $T(H)$
generated by (and consisting of) the following maps:
$$
[c;\,x,h]:\ \Phi\mapsto c\,e(h)\,\,x^*(\Phi )
\eqno(0.13)
$$
where $c\in K^*$, $x\in T(H)(K)$, $h\in H$,
and $x^*(e(h)):=h(x)\,e(h),\,h(x)$ being the value of $e(h)$
at $x$ as above.

\smallskip

In these terms, a system consisting of
a period subgroup in $T(H)(K)$ and compatible
automorphy factors becomes simply
a homomorphism, which we will call {\it a multiplier},
$$
\Cal{L}:\ B\mapsto \Cal{G}(H),\ \Cal{L}(b)=[c_b;\,x_b,h_b],
\eqno(0.14)
$$
where $B$ is a free abelian group
of the same rank as $H$ and $b\mapsto x_b$ is a bijection. (We will find convenient
in the future to weaken these restrictions).

\smallskip

The equations for formal theta functions (0.10) simply say
that $\theta$ is invariant with respect to the (image of) $B$.

\smallskip

We can also consider more general equations of the form
$$
\left(\sum_{\Gamma\in \Cal{G} (H)}c_{\Gamma }\Gamma \right)\Phi =0,\quad c_{\Gamma}\in K .
\eqno(0.15)
$$
Such equations of greater length than two
(the number of non--zero $c_\Gamma$) are satisfied by 
many classical automorphic products,
for example $e_q(t):=\prod_{n\ge 0}(1+q^{2n+1}t)$:
$$
e_q(q^2t)-(1+qt)\,e_q(t)=0.
\eqno(0.16)
$$

\smallskip

{\bf 0.5. Quantizing abelian varieties.} Let us now introduce into
the classical picture sketched above a new element,
namely, {\it a holomorphic Poisson structure $\alpha$ on $\Cal{A}$}
(and its covers) whose coefficients are constant
in a basis of invariant vector fields. 

\smallskip

In the theory of quantum linear groups $G$,  such a structure
defines a quantum deformation of all relevant algebraic data.
In particular, the function ring of $G$ is replaced
by a new non--commutative ring, with the same underlying linear
space but deformed multiplication. Comultiplication,
to the contrary, remains unchanged. In this way we get
a Hopf algebra $A_q[G]$ considered as an algebra of functions on the
appropriate quantum group.

\smallskip

We would like to play the same game with the homogeneous coordinate ring
of a Poisson abelian variety $(A,\alpha )$. In this game,
the peculiarities of multiplication and comultiplication 
look rather different from each other and from the linear case.

\smallskip

(i) {\it Multiplication and quantum tori.} Consistently with our
approach to abelian varieties, we start with quantizing the 
respective covering torus $T(H)$. This is the standard procedure:
we get a quantum torus which we denote $T(H,\alpha )$.
Over $\bold{C}$, after exponentiation (0.7)
the Poisson structure becomes an alternating bimultiplicative
pairing on $H$. Generally, over a base field $K$,
we consider a pairing $\alpha :\,H\times H\to K^*$ such that for all $h,g\in H$
$$
\alpha (h,g)=\alpha (g,h)^{-1},\quad \alpha (h_1+h_2,g)=
\alpha (h_1,g) \alpha (h_2,g).
\eqno(0.17)
$$
The space of functions of any type (algebraic,
analytic, smooth, formal ...) is topologically spanned by the set of invertible
elements, {\it formal exponents} $e(h)$. They are linearly independent over $K$ and
satisfy the deformed multiplication rule (0.6):
$$
e(h) e(h^{\prime}) =\alpha (h,h^{\prime})\,e(h+h^{\prime}).
\eqno(0.18)
$$
We may write $e_{H,\alpha}(h)$ for $e(h)$ if different
$H$ and/or $\alpha$ are discussed simultaneously.

\smallskip

As in the theory of quantum groups, $\alpha$ plays the role of 
{\it quantization parameter.} Notice, however that for nontrivial
$\alpha$,
$T(H,\alpha )$ {\it is not a quantum group}: comultiplication is lost.

\smallskip

{\it The ring of algebraic functions} $Al\,(H,\alpha )$ of the torus
$T(H,\alpha )$, by definition, is the linear space spanned
over $K$ by  $e(h),\,h\in H$ with multiplication (0.18).

\smallskip

We will also consider the two--sided $Al\,(H,\alpha )$--module
of {\it formal functions} $Af\, (H,\alpha )$ consisting
of infinite linear combinations $\sum_h a_he(h),\,a_h\in K$,
and, in the case of a complete normed field $K$ and {\it an unitary} 
quantization parameter $\alpha$ (that is, $|\alpha |=1$) the ring
of {\it analytic functions} $An\, (H,\alpha )$ consisting
of those formal functions for which $|a_h|\,\|h\|^N\to 0$
for any $N$ as $\|h\|\to \infty$, $\|h\|$ being any
Euclidean norm on $H$.

\smallskip

When $\alpha \equiv 1$, we get the usual notions of commutative geometry,
so that $T(H,1)$ is the algebraic torus previously denoted $T(H)$.
A $K$--point $x\in T(H,1)(K)$ is thus the same
as a homomorphism $H\to K^*$ given by the values
of all exponents $e_{H,1}(h)$ at $x$.

\smallskip

{\it A formal theta function} on $T(H,\alpha )$, by definition,
is an element $\Phi\in  Af\,(H,\alpha )$ invariant
with respect to an abelian group $B$ consisting
of transformations of the form
$$
[c;\,x,h_l,h_r]:\,\Phi\mapsto c\,e_{H,\alpha}(h_l)\,x^*(\Phi )\,e_{H,\alpha}(h_r)^{-1}.
\eqno(0.19)
$$
generalizing (0.10). Here $c\in K^*,\, g,h \in H$, $x\in T(H,1)(K)$,
and $x^*$ is the shift automorphism of the function
ring defined by 
$$
x^*(e_{H,\alpha}(h)):=h(x)\,e_{H,\alpha}(h),\ h(x):=e_{H,1}(x).
\eqno(0.20)
$$
Notice the all--important presence of
formal exponents at both sides of $x^*(\Phi)$ in (0.19): $l$ stands for left
and $r$ for right in the notation $h_{l,r}$. Exponents
$e_{H,\alpha}(h)$ are not central, and although 
we can transfer $e_{H,\alpha}(h_l)$ to the right or 
$e_{H,\alpha}(h_r)$ to the left, this will result
in changing $x$. 

\smallskip

In order to keep track of such changes,
it is convenient to see $\alpha$ as a homomorphism
$H\to T(H,1)(K):\,h\mapsto A_h$, defined by
$$
\forall g,h\in H,\ g(A_h)=\alpha (g,h).
\eqno(0.21)
$$
Clearly, we have $g(A_h)=h(A_{-g}),$ and any homomorphism
satisfying these conditions comes from some $\alpha$.
 
\smallskip

From (0.18) and (0.21) one sees that internal
automorphisms corresponding to formal exponents
can be alternatively seen as special shift automorphisms:
$$
e(h)\,\Phi\, e(h)^{-1}=(A_h^{-2})^*(\Phi ).
\eqno(0.22)
$$
We will call the subgroup of $T(H,1)(K)$ consisting of all
$A_h^2$ {\it the hidden period group} of the quantum
torus $T(H,\alpha ).$

\smallskip

Quantum identities for the
classical modular functions invoked in 0.7 below
clearly exhibit this double role of the parameter
$q^2$: as  the quantization parameter of the torus
and as the multiplicative period of certain automorphic
functions on this torus.

\smallskip

The quantum thetas studied in this paper
may have either hidden, or more general shift period groups.

\smallskip

Namely, imitating the classical formalism, we can define 
quantum multipliers, develop the quantum versions
of symmetry and positivity conditions
and finally, in the case of unitary $\alpha$, to
 use the multiplication of analytic
theta functions in order to introduce
a quantized version of the homogeneous coordinate
ring of an abelian variety.

\smallskip

Now, however, we do not get a ring.
We cannot say that the deformed product
of two quantum theta functions with the same
shift period groups is again such theta. In fact, what in the classical
case was a period subgroup $B$ in $T(H,1)(K)$
splits now into {\it two} period subgroups, left and right,
associated with each quantum multiplier
and which differ from each other by shifts by
hidden periods, encoded in $\alpha$.
A product of two analytic thetas is a theta,
only if the right periods of one of them coincide
with the left periods of another.

\smallskip

It turns out that
the quantum multiplication of analytic theta functions 
can be viewed as composition of morphisms 
in  a linear category, and it is this category that should  be
considered as {\it the quantum deformation of the
classical graded coordinate ring of an abelian variety.}

\smallskip

Objects of this category are {\it period
maps}: homomorphisms of abelian groups
$B\to T(H,1)(K)$. The case of embeddings of
disctrete lattices is the most important one, but 
considering arbitrary maps is imposed by our
formalism.

\smallskip

It is interesting to remark that the mirror symmetry
for abelian varieties can be formulated
exactly in terms of these period maps. Namely,
any diagram of the form:
$$
(i,i^t):\ T(H,1)(K)\leftarrow B\rightarrow T(H^t,1)(K)
$$
determines mirror dual pairs of quotients (in particular,
abelian varieties)
$(\Cal{A}:=T(H,1)/i(B), i^t)$ and
$(\Cal{B}:=T(H^t,1)/i^t(B), i)$. Here $H^t=\roman{Hom}\,(H,\bold{Z}),$
and
$i^t$ should be considered
as a version of physicists' $B$--field for $\Cal{A}$, whereas
$i$ as a $B$--field for $\Cal{B}$. 
We refer to [Ma2], Introduction
and \S 1, for further explanations.

\smallskip

It seems that our viewpoint is consistent with the remark made in [KaO]:
in the presence of a $B$--field the category of coherent sheaves
on an abelian variety must be mofified in order to fit into
the framework of homological mirror duality. We discuss only
invertible sheaves in this paper.

\medskip

(ii) {\it Comultiplication and Mumford's formalism.}
An abelian variety is, of course, a group variety. This
should provide an additional structure on its graded
function ring. This additional structure cannot, however,
be the standard Hopf comultiplication, because the pullback of an invertible
sheaf ${L}$ with respect to the
addition morphism $A\times A\to A:\,(x,y)\mapsto x+y,$
generally cannot be expressed through $pr_1^*({L})$
and $pr_2^*({L})$. D.~Mumford in [Mu] remarked that
pullbacks behave nicely with respect to the morphism
$M:\, A\times A\to A\times A:\, (x,y)\mapsto (x+y, x-y)$,
namely, $M^*({L}\boxtimes {L})\cong {L}^2\boxtimes {L}^2,$
if ${L}$ is symmetric.

\smallskip

Moreover,  $\Gamma ({L})$  carries
an irreducible representation of what we will call
{\it the small Heisenberg group} $\Cal{G} ({L})$
which in our setup can be defined as
the natural subquotient of $\Cal{G} (H,1)$ transforming
$\Gamma ({L})$ into itself.

\smallskip

Mumford developed a rich and beautiful theory showing how 
multiplication of sections of invertible sheaves and $M^*$
interacts with the action of small Heisenberg groups
and allowing him to produce explicit coordinates
on the moduli space of abelian varieties (with some
rigidity) and on the varieties themselves, and to write
the comultiplication law (or rather $M$) explicitly in terms of these
coordinates.

\smallskip

We have seen already that multiplication becomes
crucially deformed in our picture.

\smallskip

Remarkably, Mumfords morphism $M$ not only
survives quantization (recall that
already $T(H,\alpha )$ is not a quantum group), 
but also remains undeformed
in a very precise sense. We will now proceed to explaining this. 
  
\medskip

{\bf 0.6. Heisenberg  groups of quantum tori and
their rigidity.} By analogy with subsection 0.4,
we will define the Heisenberg group $\Cal{G}(H,\alpha )$
as the group generated by (an consisting of) all linear
transformations (0.19). A quantum multiplier
$\Cal{L}$ is a homomorphism of an abelian group
$B\to \Cal{G}(H,\alpha ).$ The respective space
of quantum theta functions, invariant with respect
to the image of $B$, will be denoted $\Gamma (\Cal{L}).$

\smallskip

Now, it turns out that
as any decent Heisenberg group, $\Cal{G}(H,\alpha )$
together with its basic irreducible representation $Al\,(H,\alpha )$
is a rigid object, and in particular, does not depend
on the quantization parameter. This is true in the following sense:
there is a unique ``untwisting'' isomorphism
of $\Cal{G}(H,\alpha )$ with $\Cal{G} (H,1)$ identical on the
subgroups $K^*$ and  compatible with
representations on functions which are identified
via $e_{H,\alpha}(h)\mapsto e_{H,1}(h)$ for all $h\in H$
(Proposition 1.2.1).

\smallskip

In particular, all our formal and analytic theta functions and, more generally,
solutions to the equations (0.15), are obtained by
lifting similar objects from the classical commutative world.
They are even represented by Fourier series with the same coefficients.
The small Heisenberg groups also are identified with their
commutative counterparts.
This is why Mumford's comultiplication formalism
remains unchanged, as soon as we established the existence
of the quantum version of $M$.

\smallskip

Of course, this does not mean that the theory of quantum thetas is void:
the deformed multiplication makes all the difference,
in the same way as in the theory of quantum linear groups. 

\medskip

{\bf 0.7. Geometry of non--commutative tori vs
theory of special functions.} In the spirit of A.~Connes
[Co], studying topology (smooth geometry, analytic geometry ...)
of $T(H,\alpha )$ is the same as studying modules over 
various function algebras
of $T(H,\alpha )$, their $K$--theory, connections
and Chern characters. One of the most interesting
effects is the existence of non--trivial Morita equivalences
between  $T(H,\alpha )$ with various $\alpha$, and similar
equivalences for modules with connections.

\smallskip

A recent and very intriguing motive is the appearance of non--commutative
spaces, tori in particular, as ``degenerations'' or
``extended deformations'' of usual geometric objects:
abelian and Calabi--Yau varieties, vector bundles etc.
Non--commutative tori appear also in some contemporary
versions of quantum strings and $M$--theory.

\smallskip

For a sample of these developments and further
references cf. [Ri], [RiSch],  [So], [CoDSch], [AsSch], [Wi],
[KoSo].

\smallskip

This paper is focused on another aspect of the theory,
that of special functions.

\smallskip

Physical applications of special functions are related to
statistical physics, completely integrable lattice models etc.
Roughly speaking, this happens because as we have seen
a functional equation (0.15) satisfied by the function
$F=\sum_h a_he(h)$ can be translated into a recurrent
relation between the values of the function $h\mapsto a_h$ 
or the operator function $h\mapsto a_he_{H,\alpha}(h)$ on the lattice
$H$. Even classical functions lifted in different ways
to a noncommutative torus may exhibit unexpected properties
and satisfy new functional equations.

\smallskip

One of the best known examples is the {\it $q$--exponential
function} which in our notation 
is simply $e_q(t)$ from (0.16). To see why it
deserves this name, consider the two--dimensional torus $T_q$ with the lattice $H=\bold{Z}h_1\oplus\bold{Z}h_2$
and the scalar product $\alpha (h_1,h_2)=q,$  $\alpha (h_1,h_1)=
\alpha (h_2,h_2)=1.$ Put $u=e(h_1), v=e(h_2)$ so that
$uv=q^2vu$. Then for $|q|<1$
in the ring $An\,(T_q)$ we have
$$
e_q(u)\,e_q (v)=e_q (u+v).
\eqno(0.23)
$$
Still more interesting is the formula
$$
e_q (v)\,e_q (u)=e_q (u)\,e_q (qvu)\,e_q (v) 
\eqno(0.24)
$$
derived in [FK] as the quantum version of the Rogers pentagon
identity for the classical dilogarithm. 

\smallskip

Classical theta functions
also made appearance in this context. For example, two lifts
of $\theta_q (t)$ (see (0.11)) satisfy a non--commutative identity 
$$
\theta_q (u)\,\theta_q (v)\,\theta_q (u)=
\theta_q (v)\,\theta_q (u)\,\theta_q (v).
\eqno(0.25)
$$
Notice however, that that, unlike (0.23), (0.24),  the products
(0.25) are only formal functions even for $|q|<1$. 

\smallskip

In the same vein, the function
$$
r(z,t):=\frac{\theta_q (t)}{e_q (zt)\,e_q (zt^{-1})},\ z\in K^*,
$$
satisfies the Yang--Baxter type equation
$$
r(z,u)\,r(zz^{\prime},v)\,r(z^{\prime},u)=
r(z^{\prime},v)\,r(zz^{\prime},u)\,r(z,v)
\eqno(0.26)
$$
For a proof and discussion, see [FV]. The recent preprint [FKV]
describes beautiful applications of these identities
and their generalizations to some tori of higher dimension.

\smallskip

I hope that results of this paper can be developed
in this direction and integrated with
the more geometric theory due to Connes and others.

\medskip

{\it Acknowledgement.}  General quantum thetas were first introduced
in [Ma1].
The approach presented here was developed
after the talk given by the author 
at the Workshop ``Trends in Non-commutative
Algebra and Geometry'' held at the Max Planck Institute
for Mathematics in Bonn, Sept. 25--29, 2000.
I am very grateful to the participants of this Workshop,
especially to Alain Connes, Albert Schwarz and Yan Soibelman,
discussions with whom allowed me to see more clearly
the basic features of the theory.

\bigskip

\centerline{\bf \S 1. Heisenberg groups of quantum tori}

\medskip

{\bf 1.1. Heisenberg groups and related objects.} 
In this subsection we fix
$(H,\alpha )$ as in 0.5 (i). Denote by $\widetilde{\Cal{G}} (H,\alpha )$ 
the set consisting 
of all quadruples $c\in K^*; x \in T(H,1)(K)$;
$g, h\in H,$ with the multiplication law
$$
[c^{\prime};\,x^{\prime},g^{\prime},h^{\prime}]\cdot[c;\,x,g,h]=
$$
$$
[c^{\prime}c\,g(x^{\prime})\,h(x^{\prime})^{-1}\,
\alpha (g^{\prime},g)\,\alpha (h^{\prime},h)^{-1};\,
x x^{\prime}, g+g^{\prime}, h+h^{\prime}].
\eqno(1.1)
$$
\proclaim{\quad 1.1.1. Lemma--Definition} (a) $\widetilde{\Cal{G}} (H,\alpha )$
is a group. It acts on each of the spaces $Al\,(H,\alpha )$, $An\,(H,\alpha )$,
$Af\,(H,\alpha )$ by the rule
$$
[c;\,x,g,h]:\ \Phi\mapsto
c\,e(g)\,x^*(\Phi )\,e(h)^{-1} .
\eqno(1.2)
$$

(b) Denote by $\widetilde{\Cal{G}}_0(H,\alpha )$ the subset
consisting of the elements $[1;\,A_h^2,h,h]$ for all $h\in H$
(see (0.21)).
This is a central subgroup of $\widetilde{\Cal{G}}(H,\alpha )$,
kernel of the representation (1.2).

\smallskip

(c) Put $\Cal{G}(H,\alpha )=
\widetilde{\Cal{G}}(H,\alpha )/\widetilde{\Cal{G}}_0(H,\alpha )$
and call it the large Heisenberg group of $T(H,\alpha ).$
Denote by $\Cal{G}_l(H,\alpha )$ (resp. $\Cal{G}_r(H,\alpha )$)
the subgroup of $\widetilde{\Cal{G}}(H,\alpha )$ consisting
of all elements $[c;\,x, h,0]$ (resp. $[c;\,x, 0, h]$).
Then the natural projections $\Cal{G}_{l,r}(H,\alpha )\to
\Cal{G} (H,\alpha )$ are isomorphisms.

\smallskip

(d) The map $[c;\,x,g,h]\mapsto g-h$ induces a well
defined surjective morphism $h^-:\,\Cal{G} (H,\alpha )\to H$
with kernel isomorphic to $K^*\times T(H,1)(K).$ 
\endproclaim

\smallskip

{\bf Proof.} A straightforward computation using (0.17), (0.18) shows
 that $\widetilde{\Cal{G}} (H,\alpha )$
is a group and (1.2) is a representation.

\smallskip

Clearly, $[c;\,x,g,h]$
produces the identical map, iff we have for all $k\in H$:
$$
c\,e(g)\,x^*(e(k))\,e(h)^{-1}=e(k).
$$
This means that $g=h$ and  for all $k$,
$c\,k(x)\,\alpha (h,k)^2\,=1,$
in other words,
$c=1,\ k(x)=\alpha (k,h )^2 .$
Comparing this to (0.21), one sees that $x=A_h^2$. 

\smallskip

Using this description, one easily checks that
$\widetilde{\Cal{G}}(H,\alpha )$ is a semidirect product
of either of the groups $\Cal{G}_{l,r}(H,\alpha )$ 
and of the central subgroup $\widetilde{\Cal{G}}_0(H,\alpha )$, which finishes the proof. 

\smallskip

Many properties of $\Cal{G} (H,\alpha )$ can be conveniently checked
by using its identification with, say, $\Cal{G}_l(H,\alpha )$.

\smallskip

However, in order to treat the multiplication of theta functions below,
we will have to use both versions $\Cal{G}_l$ and $\Cal{G}_l$
on the equal footing. To this end we will collect several
formulas for future use. Call the left (resp. the right)
representative of $[c;x,g,h]$ its natural projection of the form
$[c_l;x_l,h_l,0]$ (resp. $[c_r;x_r,0, h_r]$). We have then
$$
[c_l;x_l,h_l,0]=[c\,\alpha (h,g)\varepsilon (h);xA_h^{-2},g-h,0],
\eqno(1.3)
$$
$$
[c_r;\,x_r,0,h_r]=[c\,\alpha (h,g)\varepsilon (g);xA_g^{-2},0,h-g].
\eqno(1.4)
$$
Here $\varepsilon (h):=\alpha (h,h)$; this is a character of $H$
with values in $\{\pm 1\}$ which in [Ma2] was called
{\it the characteristic} of the quantization parameter $\alpha$.

\smallskip

To put this differently, two elements $[c_l;x_l,h_l,0],\,[c_;x_r,0,h_r]$
have the same image in $\Cal{G} (H,\alpha )$ iff
$$
h_l=-h_r,\ x_r=x_l A_{h_l}^{-2}=x_l A_{h_r}^{2},\
c_r=c_l\varepsilon (h_l)=c_l\varepsilon (h_r).
\eqno(1.5)
$$
\smallskip

{\bf 1.2. Rigidity.} As we mentioned in the Introduction, 
our Heisenberg group
together with its basic irreducible representation
is a rigid object independent of $\alpha$. More precisely,
let us refer to the quantization parameter explicitly
by adding  an appropriate subscript to our notation.
Then we have:

\smallskip

\proclaim{\quad 1.2.1. Proposition} (a) The map written
in terms of left representatives
$$
u_{1,\alpha}:\,\Cal{G}(H,1)\to \Cal{G}(H,\alpha ):\
[c;\,x, h,0]_1\mapsto [c;\,xA_h, h,0]_{\alpha}
\eqno(1.6) 
$$
is a group isomorphism.

\smallskip

(b) The linear map (denoted by the same symbol) 
$$
u_{1,\alpha}:\,Al\,(H,1)\to Al\,(H,\alpha):\
e_{H,1}(h)\mapsto e_{H,\alpha}(h)
\eqno(1.7)
$$
is an isomorphism of the representations extending (1.6),
and similarly for other spaces of functions.

\smallskip

(c) $u_{1,\alpha}$ is the unique isomorphism satisfying (b).
\endproclaim

\smallskip

{\bf Proof.} In fact, one easily checks that
the map (1.6) is bijective and compatible with products.

\smallskip

Moreover,
$$
[c;\,x, h,0]_1\,(e_{H,1}(g))= c\,e_{H,1}(h)\,x^*(e_{H,1}(g))=
c\,g(x)\,e_{H,1}(h+g),
$$
$$
[c;\,xA_h, h,0]_{\alpha}\,(e_{H,\alpha }(g))= c\,e_{H,\alpha}(h)\,(xA_h)^*(e_{H,\alpha}(g))=
c\,g(xA_h)\,\alpha (h,g)\, e_{H,\alpha}(h+g),
$$
and the coefficients at the right hand sides  coincide
in view of (0.21). Finally, unicity follows from the fact that
our representation is faithful.

This proves the Proposition.

\medskip

More generally, we can define a coherent system of
twisting isomorphisms
$$
u_{\alpha ,\beta}:=u_{1,\beta}\circ u_{1,\alpha}^{-1}:\,
\Cal{G}_l(H,\alpha )\to \Cal{G}_l(H,\beta ):\
[c;\,x, h,0]_{\alpha}\mapsto [c;\,xA_h^{-1}B_h, h,0]_{\beta},
\eqno(1.8) 
$$
where the points $B_h$ are determined by the relations (0.21)
with $\beta$ replacing $\alpha$.

\medskip

Of course, twisting isomorphisms generally are not compatible with
multiplication. We will now introduce and study the structure
on the Heisenberg group which reflects multiplication of
functions and which actually depends on $\alpha$.

\medskip

{\bf 1.3. Partial composition in $\Cal{G} (H,\alpha ).$}
Let $\Gamma^{\prime\prime}, \Gamma^{\prime}\in \Cal{G} (H,\alpha ).$
Call these elements {\it composable} (in this order)
iff $x^{\prime\prime}_r=x^{\prime}_l$ where $x^{\prime\prime}_r$
(resp.  $x^{\prime}_l$) is taken from the right (resp. left)
representative of $\Gamma^{\prime\prime}$ (resp. $\Gamma^{\prime}$):
$$
\Gamma^{\prime\prime} \equiv [c_r^{\prime\prime};\,x^{\prime\prime}_r,0,h^{\prime\prime}_r]\
\roman{mod}\,\widetilde{\Cal{G}}_0 (H,\alpha ),
$$
$$
\Gamma^{\prime} \equiv [c_l^{\prime};\,x^{\prime}_l,h^{\prime}_l,0]\
\roman{mod}\,\widetilde{\Cal{G}}_0 (H,\alpha ).
$$

\smallskip

If this condition is satisfied, define the composition
$$
\Gamma^{\prime\prime}\circ \Gamma^{\prime}:=
[c_r^{\prime\prime};\,1,0,h^{\prime\prime}_r]\cdot
[c^{\prime}_l;\,x^{\prime}_l,h^{\prime}_l,0] 
\,\roman{mod}\,\widetilde{\Cal{G}}_0 (H,\alpha ).
\eqno(1.9)
$$
Calculating with the help of (1.1), (1.3), (1.4), we get
$$
\Gamma^{\prime\prime}\circ \Gamma^{\prime}\equiv
[c^{\prime\prime}_rc^{\prime}_l\,\alpha (h^{\prime}_l,
h^{\prime\prime}_l)\,\varepsilon (h^{\prime\prime}_l);\,x^{\prime}_lA^2_{h_l^{\prime\prime}},
h^{\prime}_l+h^{\prime\prime}_l,0]\equiv
$$
$$
[c^{\prime\prime}_rc^{\prime}_l\,\alpha (h^{\prime}_l,
h^{\prime\prime}_l)\,\varepsilon (h^{\prime}_l);\,x^{\prime\prime}_rA^2_{h_r^{\prime}},0,
h^{\prime}_r+h^{\prime\prime}_r]\
\roman{mod}\,\widetilde{\Cal{G}}_0 (H,\alpha ) .
\eqno(1.10)
$$
\smallskip

\proclaim{\quad 1.3.1. Proposition} (a) Assume that the pairs
$(\Gamma^{\prime\prime\prime}, \Gamma^{\prime\prime})$ and
$(\Gamma^{\prime\prime}, \Gamma^{\prime})$
are composable. Then the pairs
$(\Gamma^{\prime\prime\prime},\Gamma^{\prime\prime}\circ\Gamma^{\prime})$
and 
$(\Gamma^{\prime\prime\prime}\circ\Gamma^{\prime\prime},\Gamma^{\prime})$
are composable as well, and the triple composition
is associative.

\smallskip

(b) We have for any $\Phi,\Psi\in Al\, (H,\alpha )$
(or $An\, (H,\alpha )$)
$$
\Gamma^{\prime} (\Phi )\,\Gamma^{\prime\prime} (\Psi )=
(\Gamma^{\prime\prime}\circ\Gamma^{\prime})\, (\Phi\Psi ).
\eqno(1.11)
$$
\endproclaim

\smallskip

{\bf Proof.} We start with checking compatibilities
needed to form two triple products. Given equalities
$x^{\prime\prime\prime}_r=x^{\prime\prime}_l$ and
$x^{\prime\prime}_r=x^{\prime}_l$, we obtain,  using (1.10):
$$
x_l^{\prime\prime\circ\prime}=x^{\prime}_lA^2_{h^{\prime\prime}_l}=
x^{\prime\prime}_rA^2_{h^{\prime\prime}_l}
$$
which in view of (1.5) coincides with $x_l^{\prime\prime}$ 
and therefore with $x_r^{\prime\prime\prime}$. Here
$x_l^{\prime\prime\circ\prime}$ denotes the $x$--component
of the left representative of $\Gamma^{\prime\prime}\circ \Gamma^{\prime}$.

\smallskip

Similarly, $x_r^{\prime\prime\prime\circ\prime\prime}=
x^{\prime}_l$ because both sides are equal to
$x^{\prime\prime}_lA^{-2}_{h^{\prime\prime}_l}=x_r^{\prime\prime}.$

\smallskip

Now we verify (1.11):
$$
\Gamma^{\prime} (\Phi )\,\Gamma^{\prime\prime} (\Psi )=
c^{\prime}_le(h^{\prime}_l)\,x^{\prime *}_l(\Phi )\,
c^{\prime\prime}_r x^{\prime\prime *}_r(\Psi )\,
e(h^{\prime\prime}_r)^{-1}=
$$
$$
c^{\prime}_lc^{\prime\prime}_re(h^{\prime}_l)\,x^{\prime *}_l(\Phi\Psi )
e(h^{\prime\prime}_r)^{-1} =
$$
$$
[c_r^{\prime\prime};\,1,0,h^{\prime\prime}_r]\cdot
[c^{\prime}_l;\,x^{\prime}_l,h^{\prime}_l,0]\, (\Phi\Psi )=
(\Gamma^{\prime\prime}\circ \Gamma^{\prime})\, (\Phi\Psi )
$$
in view of (1.9).

\smallskip

Finally, applying (1.11) several times, we see that
$$
(\Gamma^{\prime\prime\prime}\circ (\Gamma^{\prime\prime}\circ\Gamma^{\prime}))
\,(\Phi \Psi \Xi )=
\Gamma^{\prime}(\Phi )\,\Gamma^{\prime\prime}(\Psi )\,
\Gamma^{\prime\prime\prime}(\Xi )=
((\Gamma^{\prime\prime\prime}\circ \Gamma^{\prime\prime})\circ\Gamma^{\prime})
\,(\Phi\Psi\Xi )
$$
for any $\Phi ,\Psi ,\Xi$, which proves associativity.

\medskip

We will show now that
$\Cal{G}(H,\alpha )$ with its partial multiplication can be treated
as the set of morphisms of an amusing category $\Cal{C}(H,\alpha )$.
By definition, its objects are points $\xi\in T(H,1)(K)$ and
morphisms are 
$$
\roman{Hom}\,(\xi ,\eta ):=\{\Gamma\in\Cal{G} (H,\alpha )\,|\,
x_r(\Gamma )=\xi,\, x_l(\Gamma )=\eta\}.
\eqno(1.12)
$$
For composition we take $\circ$. Proposition 1.3.1 shows that
this category is well defined. Identity morphisms are 
$[1;\,\xi,0,0].$

\medskip
\proclaim{\quad 1.3.2. Proposition} $\Cal{C}(H,\alpha )$
is a groupoid with the set of isomorphism
classes of objects $T(H,1)(K)/P_{\alpha}$,
where $P_{\alpha}:=\{A_h^2\,|\,h\in H\}$ is the group
of hidden periods. The automorphism
group of any object is isomorphic to $K^*\times K_0$ where
$K_0$ is the kernel of the form $\alpha^2$.
\endproclaim

\smallskip

{\bf Proof.} If there is a morphism $\Gamma :\,\xi\to \eta$,
represented by $[c_r;\,x_r,0,h_r]$ and $[c_l;\,x_l,h_l,0]$,
we must have $x_r=\xi ,x_l=\eta$.  In view of (1.5) this
is possible precisely when $\xi \equiv \eta\,\roman{mod}\, P_{\alpha}.$
Moreover, we then have $\xi =\eta A_{h_l}^{-2}=\eta A_{h_r}^{2}$.

\smallskip

In this case we van solve the equation $\Gamma^{\prime}\circ\Gamma =
[1;\xi ,0,0]$ for the left inverse to $\Gamma$: namely
$\Gamma^{\prime}=[c_l^{-1};\,\eta ,-h_l,0]$, and similarly for the right inverse.

\smallskip

Finally, the (left representatives of the) automorphism
group of $\xi$ are $[c;\,\xi, h,0]$, $h\in \roman{Ker}\,\alpha^2,$ with the multiplication law
$$
[c^{\prime};\,\xi, h^{\prime},0]\circ [c;\,\xi, h,0]=
[c^{\prime}c;\,\xi, h^{\prime}+h,0].
$$

\medskip

{\bf 1.3.3. Double--sided representatives.} Elements of the form
$[c;\,1,g,h]$ form a subgroup $\Cal{G}_d (H,\alpha )$
in $\widetilde{\Cal{G}} (H,\alpha ).$ 
Assume that
$\alpha^2$ is non--degenerate. Then the natural morphism
$\Cal{G}_d (H,\alpha )\to {\Cal{G}} (H,\alpha )$ is an embedding
(see Lemma 1.1.1 (b)). We will identify in this case
$\Cal{G}_d (H,\alpha )$ with its image, and call its elements
double--sided representatives.

\smallskip

Two such representatives $[c^{\prime\prime},1,h^{\prime\prime}_l,
h^{\prime\prime}_r]$, $[c^{\prime},1,h^{\prime}_l,
h^{\prime}_r]$ are composable, iff $h_r^{\prime}=h_l^{\prime\prime}$,
and in this case their composition is double--sided as well:
$$
[c^{\prime\prime},1,h^{\prime\prime}_l,
h^{\prime\prime}_r]\circ [c^{\prime},1,h^{\prime}_l,
h^{\prime}_r] = [c^{\prime\prime}c^{\prime};\,h^{\prime}_l,
h^{\prime\prime}_r].
$$
\smallskip 
We leave the straightforward check as an exercise.

\medskip

{\bf 1.4. Functorial properties of the Heisenberg group.}
In this paper as in [Ma2] the category of quantum tori
is defined as opposite to the category of their
rings of {\it algebraic} functions. This means
that a morphism $F:\,T(H,\alpha )\to T(H^{\prime},\alpha^{\prime})$,
is given by the contravariant $K$--algebra homomorphism
$F^*:\,Al\,(H^{\prime},\alpha^{\prime})\to Al\,(H,\alpha )$.
The following (easy) result is proved in [Ma2]:

\medskip

\proclaim{\quad 1.4.1. Proposition} a) The set of invertible elements
of $Al\,(H,\alpha )$ is $\{ a\,e(h)\,|\,a\in K^*, h\in H\}.$
If  $F:\,T(H_2,\alpha_2 )\to T(H_1,\alpha_1)$ is a morphism of non--commutative tori, then the induced map $f=[F]:\, H_1\to H_2$
determined by $F^*(e_{H_1,\alpha_1}(h))=a_he_{H_2,\alpha_2}(f(h)),\, a_h\in K^*,$ 
is additive and compatible with the squares of the quantization
forms:
$$
\alpha_2^2(f(h),f(g))=\alpha_1^2(h,g).
\eqno(1.13)
$$
\smallskip

b) The set of all morphisms $F$ with a fixed $f=[F]$ is either empty, or
has a natural structure of the principal homogeneous space
over the group $T(H_1,1)(K)=\roman{Hom}\,(H_1,K^*).$
In particular, if the characteristic of $f$ is 1,
then $F^*:\,e(h)\mapsto e(f(h))$ extends to a uniquely
defined morphism of rings of algebraic functions.

\smallskip

c) Any morphism $F^*$ extends to $Af$ by $F^*(\sum a_h e(h))=
\sum a_hF^*(e(h)).$ If $K$ is normed and $\alpha$ unitary,
then this extension maps analytic functions to analytic.
\endproclaim

The bilinear form 
$$
\varepsilon_f (h,g):= \alpha_1(h,g)\alpha_2^{-1}(f(h),f(g))
$$
with values in $\{\pm 1\}$ is called the characteristic of $f$
and of $F$. This form is trivial iff
$a_{g+h}=a_ga_h$ for all $g,h$.

\medskip

Shifts (0.20) of $T(H,\alpha )$ by points of $T(H,1)(K)$ are endomorphisms of $T(H,\alpha )$. Besides them,
we will use the {\it multiplication by $n$}
$$
[n]:\, T(H,\alpha )\to T(H,\alpha^{n^2}),\quad
[n]^*(e_{H,\alpha^{n^2}}(h))=e_{H,\alpha}(nh).
\eqno(1.14)
$$
and {\it the Mumford morphism}
$$
M:\,T(H\oplus H,\alpha\oplus\alpha )\to T(H\oplus H,\alpha^2\oplus\alpha^2),
$$
$$
M^*(e_{H\oplus H,\alpha^2\oplus\alpha^2}(h,g))=
e_{H\oplus H,\alpha\oplus\alpha }(h+g,h-g).
\eqno(1.15)
$$
It is well defined, because
$$
(\alpha\oplus\alpha )[(h+g,h-g),(h^{\prime}+g^{\prime}, h^{\prime}-g^{\prime})]=
$$
$$
\alpha (h+g, h^{\prime}+g^{\prime}) \alpha (h-g, h^{\prime}-g^{\prime})
=\alpha^2(h,h^{\prime}) \alpha^2(g,g^{\prime}).
$$
\smallskip

Consider now a morphism of tori $F:\,T(H_2,\alpha_2 )\to T(H_1,\alpha_1)$
as in Proposition 1.4.1 such that $F^*(e(h))=a_he(f(h)),\, a_h\in K^*,$
for a homomorphism of character groups $f=[F]:\, H_1\to H_2$.
Assume that $\roman{Ker}\,f = \{0\}$ so that $F^*$ is an embedding.
Denote by $\Cal{G} (F)$ the following subquotient
of $\Cal{G} (H_2,\alpha_2)$: take the full subgroup
stabilizing the subspace $F^*(Al\,(H_1,\alpha_1))$
and mode out the kernel of this representation.

\medskip 

\proclaim{\quad 1.4.2. Proposition} Assume that $F^*$
satisfies the following condition: $h\mapsto a_h$
is a group homomorphism. Then there exists a well defined
isomorphism $\Cal{G} (F)\to \Cal{G} (H_1, \alpha_1)$ compatible
with the natural representations of these groups on $Al\,(H_1,\alpha_1)$.
\endproclaim

\smallskip

{\bf Proof.} Let us work with left representatives.

\smallskip

Clearly $[c;\,x, h,0]_{\alpha_2}$ stabilizes the subspace
$F^*(Al\,(H_1,\alpha_1))$ iff $h\in \roman{Im}\,f$.
Such an element acts identically on this subspace
iff $c=1$, $h=0$, and $f^*(g)(x)=1$ for all $g\in H_1$, in other
words, $x$ belongs to the kernel $\Cal{K}$ of the
homomorphism $\varphi :\,T(H_2,1)(K)\to T(H_1,1)(K)$ induced by $f$.

\smallskip

For any $g,h\in H_1$ we then have
$$
[\,ca_g;\, x, f(g),0]_{\alpha_2}\, F^*(e(h))= F^*([\,c;\,\varphi (x),g,0]_{\alpha_1} \,e(h))
$$
(here the multiplicativity of $a_h$ is used). One now easily
sees that the map
$$
\Cal{G} (F)\to \Cal{G} (H_1, \alpha_1)\,:\
[\,ca_g;\, x\,\roman{mod}\,\Cal{K}, f(g),0\,]_{\alpha_2}\mapsto
[\,c;\,\varphi (x),g,0\,]_{\alpha_1}
$$
is a well defined isomorphism which we looked for.

\medskip

{\bf 1.4.3. Homomorphisms $\psi_{n,d}$.} Let $d,n$ be two
integers, $d/n$. Consider the map
$$
\psi_{d,n}:\,\widetilde{\Cal{G}} (H,\alpha^d)\to 
\widetilde{\Cal{G}} (H,\alpha ):\ 
[c;\,x,h_l,h_r]_{\alpha^d}\mapsto [c^{n^2/d};\,x^{n/d},nh_l,nh_r]_{\alpha}.
\eqno(1.16)
$$
A straightforward check shows that it is a homomorphism
mapping $\widetilde{\Cal{G}}_0 (H,\alpha^d)$ into
$\widetilde{\Cal{G}}_0 (H,\alpha )$. Hence it induces
a homomorphism of the large Heisenberg groups which we denote
by the same symbol:
$$
\psi_{d,n}:\,\Cal{G} (H,\alpha^d)\to 
\Cal{G} (H,\alpha ).
$$

\newpage

\centerline{\bf \S 2. Multipliers and theta functions}

\medskip

\proclaim{\quad 2.1. Definition} (i) A (formal) theta multiplier 
for $T(H,\alpha )$ is a homomorphism $\Cal{L}:\,
B\to \Cal{G} (H,\alpha )$ where $B$ is an abelian group.

\smallskip

(ii) A theta function with multiplier  $\Cal{L}$ is a formal
function on $T(H,\alpha )$ invariant with respect to the
action of (the image of) $B$.

\smallskip

(iii) $\Gamma (\Cal{L})$ is the linear space of theta functions with multiplier  $\Cal{L}$. 
\endproclaim

\smallskip

This notation
is supposed to remind the case of usual abelian
varieties where we deal with invertible sheaves and their sections.

\medskip

Only the image of $B$ determines the space of theta--functions, but
using arbitrary homomorphisms gives more flexibility.

\medskip

Image of any element $b\in B$ has left and right representatives
in $\Cal{G}_{l,r}$ which we will denote respectively
$$
[\,c_{l,b};\,x_{l,b},h_{l,b},0],\  [\,c_{r,b};\,x_{r,b},0,h_{r,b}].
$$
As in (1.5), their coincidence modulo $\Cal{G}_0 (H,\alpha )$ means that
for all $b\in B$
$$
-h_{r,b}=h_{l,b}=h^-_b,\ x_{r,b}=x_{l,b}\,A^{-2}_{h_{l,b}},\
c_{r,b} =\varepsilon (h_b^-)\,c_{l,b}.
\eqno(2.1)
$$

\smallskip

\proclaim{\quad 2.2. Lemma} Let $B$ be an abstract
abelian group. Consider two families of elements 
$[\,c_{l,b};\,x_{l,b},h_{l,b},0]$ and
$[\,c_{r,b};\,x_{r,b},0,h_{r,b}]$ 
in $\Cal{G}_{l,r}(H,\alpha )$ respectively, indexed by $B$
and satisfying (2.1). Then these families 
constitute a theta multiplier as above, iff
the following conditions are satisfied:

\smallskip

(i) The maps $B\to H:\, b\mapsto h_{l,b}$ and  and $B\to T(H,1)(K):\,
b\mapsto x_{l,b}$ are group homomorphisms, and the same for
$h_r$, $x_r$.

\smallskip

(ii) The map 
$$
B\times B\to K^*:\ (b_1,b_2)\mapsto \langle b_1,b_2\rangle
:=h^-_{b_2}(x_{l,b_1})\,\alpha (h^-_{b_1},h^-_{b_2})
\eqno(2.2)
$$
is a symmetric bimultiplicative pairing. 

\smallskip

(iii) We have identically
$$
\frac{c_{l,b_1+b_2}}{c_{l,b_1}c_{l,b_2}}=\langle b_1,b_2\rangle 
\eqno(2.3)
$$
and the same for $c_{r,b}$.
\endproclaim

\smallskip

The proof is straightforward.

\medskip

If a symmetric bimultiplicative pairing $(b_1,b_2)$
is chosen such that $(b_1,b_2)^2= \langle b_1,b_2\rangle$,
then any solution of (2.3) can be uniquely represented
in the form
$$
c_{l,b}=\psi_l(b)\, (b,b), 
\eqno(2.4)
$$
where $\psi_l :\, B\to K^*$ is a character. We have $\psi_r(b)=
\varepsilon (h^-_b) \psi_l(b).$

Such choice of a square root becomes always possible
after a finite extension of $K$. We will also use the 
explicit representation of coefficients in the form (2.4),
because it matches the classical notation.

\medskip
\proclaim{\quad 2.3. Definition} (i) Quadruples
$(\psi_l, (\,,), x_l,h_l)$ (resp. $(\psi_r, (\,,), x_r,h_r)$)
are called the left (resp. right) automorphy factors
of the multiplier $\Cal{L}$.

\smallskip

(ii) The maps $x_{l,r}:\, B\to T(H,1)(K)$ and their images are called
the left (resp. the right) period groups of $\Cal{L}$.
\endproclaim

\medskip

The abstract definition of multipliers 2.1 makes 
transparent their functorial properties, whereas
the description 2.2 provides the link with the
classical constructions in the theory of theta functions.
In particular, (2.2) for $\alpha \equiv 1$ becomes
the Riemann symmetry condition so that generally
(2.2) represents its quantum deformation. 

\smallskip

Twisting isomorphisms (1.8) identify the sets
of theta multipliers for various $\alpha$'s
and their actions on the spaces of functions. 
What remains specific for each $\alpha$,
is the interaction of a multiplier with the
deformed multiplication. The period groups also
change under twisting isomorphisms. 

\medskip

{\bf 2.4. Analytic theta functions and Riemann's
symmetry conditions.} By definition, a theta
function with multiplier $\Cal{L}$ must satisfy the functional equations
$$
\psi_l (b)\,(b,b)\,e(h_{l,b})\,(x_{l,b})^*(\theta ) =\theta =
\psi_r (b)\,(b,b)\,(x_{r,b})^*(\theta )\,e(h_{r,b})^{-1}
\eqno(2.5)
$$
for all $b\in B.$ 

\medskip

\proclaim{\quad 2.4.1. Theorem} (a)  We have 
$$
\roman{dim}\,\Gamma (\Cal{L}) \le [H:h^-(B)].
\eqno(2.6)
$$
Equality holds if $x_{l,b}, c_{l,b}$ depend only on 
$h^-_b$, for example, if 
$h^-:\,B\to H$ is an embedding.

\smallskip

(b) Let the last condition be satisfied.
Assume that $K$ is a normed field.
Then all theta functions
with multiplier $\Cal{L}$ are analytic if $[H:h^-(B)] <\infty$
and $\roman{log}\,|\,\langle b,b\rangle |$
is a negatively defined
quadratic form on $B$. 

In this case
$\roman{rk}\,B = \roman{rk}\,H$, and
$B$ is discrete in $T(H,1)(K)$.

\smallskip

Such multipliers will be called ample.

\endproclaim 

\smallskip

{\bf Proof.} Let $\theta =\sum_{h\in H} a_he(h),\,a_h\in K,$ $b\in B.$
We have
$$
(x_{l,b})^*(\theta )=\sum_{h\in H}a_hh(x_{l,b})e(h)=\sum_{h\in H}
a_{h-h^-_{b}}\,h_b^-(x_{l,b})^{-1}\,h(x_{l,b})\,e(h-h^-_{b}),
$$
so that the left hand side of (2.5) becomes
$$
c_{l,b}\,e(h_{b}^-) \sum_{h\in H}
a_{h-h^-_{b}}\,h_b^-(x_{l,b})^{-1}\,h(x_{l,b})\,e(h-h^-_{b})=
$$
$$
\sum_{h\in H} c_{l,b}\, a_{h-h^-_{b}}\,h_b^-(x_{l,b})^{-1}\,h(x_{l,b})\,
\alpha (h_b^-,h)\,\varepsilon (h_b^-)\, e(h).
$$
Replacing here $h_b^-(x_{l,b})^{-1}$ by 
$\langle b,b\rangle\,\varepsilon (h^-_b)$
(see (2.2)), we see that (2.5) is equivalent to the
following equations for coefficients:
$$
a_{h-h^-_b}=a_h\, c_{l,b}^{-1}\,\langle b,b\rangle\, 
h(x_{l,b})^{-1}\,
\alpha (h,h_b^-) 
\eqno(2.7)
$$
for all $h\in H$ and $b\in B.$ Therefore if values $a_h$
for some system of representatives of $H/h^-(B)$
are fixed, $\theta$ is defined uniquely, if it exists
at all.

\smallskip

Such a choice can be done arbitrarily, if $h^-$ is injective; 
otherwise
(2.7) might be overdetermined. This proves the first statement of the theorem.
It also shows that if $\Gamma (\Cal{L})$ is not finite dimensional,
it necessarily contains non--analytic functions.
\smallskip

Assume now that $\Gamma (\Cal{L})$ is finite dimensional.
Then on each coset $h-h^-(B)$ we have, in the notation (2.4),
$$
\roman{log}\,|a_{h-h^-_b}| =
\roman{log}\,|a_h| +  \roman{log}\,|(b,b)|
- \roman{log}\,|\psi_l(b)\, h(x_{l,b})\,\alpha (h_{b}^+,h)|.
$$
The second summand in the right hand side is quadratic in $b$
whereas the third is linear. Hence analyticity follows from
the negative definiteness of $(b,b)$ or $\langle b,b \rangle$
which is the same.

\medskip

{\bf 2.5. Operations on multipliers.} {\it (i) Powers.} We can compose
any multiplier $\Cal{L}:\, B\to \Cal{G} (H,\alpha^d)$
with the homomorphism $\psi_{d,n}:\, \Cal{G} (H,\alpha^d)\to
\Cal{G} (H,\alpha)$ defined in 1.4.3. The resulting composition
will be a multiplier as well.

\smallskip

Especially important is the case $d=n$. We will denote
$\psi_{n,n}\circ\Cal{L}$ by $\Cal{L}^n$. From (1.16)
we see that
$$
\Cal{L}^n(b)=[c_{l,b}^n;\,x_{l,b}, nh_{l,b},0].
\eqno(2.8)
$$
In particular, the left and the right period groups for
$\Cal{L}^n$ are the same as for $\Cal{L}$.

\smallskip

If $\alpha\equiv 1$, $\Cal{L}\mapsto\Cal{L}^n$ simply corresponds 
to the $n$--th
power of invertible sheaves, in particular,
product of theta functions induces a map $\Gamma (\Cal{L})^{\otimes n}\to
\Gamma (\Cal{L}^n).$ This cannot be true in general, because
$\Cal{L}$ and $\Cal{L}^n$ correspond to different quantization parameters.
The only polylinear map of this kind I can think of
can be obtained by first untwisting $\alpha^n$, multiplying
in the commutative ring, and then again twisting to $\alpha$.

\smallskip

{\it (ii) External tensor products.} Consider two
multipliers 
$\Cal{L}^{\prime}:\,B^{\prime}\to
T(H^{\prime},\alpha^{\prime} )$, and $\Cal{L}^{\prime\prime}:\,
B^{\prime\prime}\to
T(H^{\prime\prime},\alpha^{\prime\prime} )$. We will use  their
left representatives denoted as in 2.1. 

\smallskip

Their external
tensor product $\Cal{L}=\Cal{L}^{\prime}\boxtimes \Cal{L}^{\prime\prime}$
is a multiplier  
$$
\Cal{L}^{\prime}\boxtimes \Cal{L}^{\prime\prime}:\
B^{\prime}\oplus B^{\prime\prime}: \to \Cal{G}(H^{\prime}\oplus H^{\prime\prime}, \alpha^{\prime}\oplus\alpha^{\prime\prime})
$$
which is defined by the following left representatives:
$$
\Cal{L} (b_1,b_2)=
[c_{l,b_1}^{\prime} c_{l,b_2}^{\prime\prime};\,
(x_{l,b_1}^{\prime}, x_{l,b_2}^{\prime\prime}),
(h_{l,b_1}^{\prime}, h_{l,b_2}^{\prime\prime}),0].
\eqno(2.9)
$$
\smallskip

{\it (iii) Pullbacks.} Consider
a morphism of non--commutative tori $F:\,T(H_2,\alpha_2 )\to T(H_1,\alpha_1)$
such that $F^*$ is an embedding and the map $h\mapsto a_h$
is a homomorphism, as
in Proposition 1.4.2. Let $\Cal{K}$ be the group defined
in the proof of that Proposition. Denote by
$\widetilde{\Cal{G}}(F)$ the subgroup of $\Cal{G}(H_2,\alpha_2)$
stabilizing $F^*(Al\,(H_1,\alpha_1))$ so that we have
the exact sequence
$$
1\to \Cal{K}\to \widetilde{\Cal{G}}(F)\to \Cal{G}(F) \to 1.
$$
A multiplier $\Cal{L}:\,B\to \Cal{G} (H_1,\alpha_1)$ 
can be identified with a homomorphism
$\Cal{L}:\,B\to \Cal{G} (F).$ Any lift of
this homomorphism to a homomorphism 
$F^*(\Cal{L}):\,B\to \widetilde{\Cal{G}} (F)$ 
is a multiplier for $T(H_2,\alpha_2).$ Such lifts certainly
exist if $\Cal{K}$ is finite and the subgroup of shifts
$\{x_{l,b}\}$ of $\Cal{L}$ is free.

\smallskip

For any lift, any function $\Phi$ on $T(H_1,\alpha_1)$ and any $b\in 
B$
we will have
$$
F^*(\Cal{L})(b) (F^*(\Phi ))=F^*(\Cal{L}(b)\Phi ).
\eqno(2.10)
$$
As a corollary, we get
$$
F^*(\Gamma (\Cal{L}))\subset \Gamma (F^*(\Cal{L})).
\eqno(2.11)
$$
\smallskip

This construction is applicable, in particular, to shifts, 
multiplications $[n]$  and Mumford's morphisms (see (1.14) and (1.15)).

\smallskip

Concretely, if as usual we denote by $[c_{l,b};\,x_{l,b},h_{l,b},0]$
the left representative of $\Cal{L}(b)$, then
$$
F^*(\Cal{L})(b)=[c_{l,b}\,a_{h_{l,b}};\, x^{\prime}_{l,b},f(h_{l,b}),
\eqno(2.12)
$$
where as above, $F^*(e(h))=a_h e(f(h)),$ $\varphi$ is the morphism
induced on commutative tori, and $x^{\prime}_{l,b}\in \varphi^{-1}(x_{l,b}).$

\smallskip

Applying this to shift morphisms and $[n]$, we get
$$
y^*(\Cal{L})(b)=[c_{l,b}h_{l,b}(y);\,x_{l,b}y^{-1},h_{l,b},0],
\eqno(2.13)
$$
$$
[n](\Cal{L})(b)=[c_{l,b};\,x_{l,b}^{1/n},nh_{l,b},0]_{\alpha^{n^2}}.
\eqno(2.14) 
$$
(Note the change of the quantization parameter in the last formula).

\smallskip

The multiplier $\Cal{L}$ is called {\it symmetric}, if 
$[-1]\Cal{L}(-b)\equiv \Cal{L}(b).$ In view of (2.14), this is equivalent
to $c_{l,b}\equiv c_{l,-b}$, or, in terms of the automorphy
factors (2.4), $\psi (b)\in \{\pm 1\}.$

\medskip

{\bf 2.6. Twisting isomorphisms.} The existence and unicity
of the twisting maps (1.6), (1.7), (1.8) shows that most of 
the structures described in this section for different
$\alpha$ are connected
by the canonical isomorphisms. More precisely, if $\Cal{L}$
is a (formal) theta multiplier for $T(H,\alpha )$, then
$u_{\alpha ,\beta}\circ \Cal{L}$ is a theta multiplier
for $T(H,\beta )$. Moreover, $\Gamma (\Cal{L})$
turns into $\Gamma (u_{\alpha ,\beta}\circ \Cal{L})$ by
simply replacing $e_{H,\alpha}(h)$ with $e_{H,\beta}(h)$.

\smallskip

The behaviour of multipliers with respect to 
morphisms of quantum tori becomes slightly more complex, if
shifts are involved. However, any commutative diagram
involving only morphisms which transform formal
exponents to formal exponents makes sense and remains commutative
after untwisting all quantization parameters to
$1$. Introducing theta functions into this picture,
we can coherently untwist them as well.

\smallskip

The most important application of this remark
concerns Mumford's morphisms (1.15), to which
I will return in \S 4. It follows
that the ``comultiplication'' of quantum theta functions
remains essentially the same as in the classical case.

\smallskip

We now pass to the more interesting effects
related to multiplication.

\bigskip

\centerline{\bf \S 3. Composition of multipliers and multiplication of thetas}

\medskip

{\bf 3.1. Composition of multipliers.} Consider
two multipliers $\Cal{L}^{\prime\prime}\,,\Cal{L}^{\prime} :\
B\to \Cal{G} (H,\alpha ).$ Call them {\it composable},
if for every $b\in B$,  $\Cal{L}^{\prime\prime}(b)$
and $\Cal{L}^{\prime}(b)$ are composable
in the sense of 1.3. 

\smallskip

When $\Cal{L}^{\prime\prime}$ and $\Cal{L}^{\prime}$ are composable,
we can define their pointwise composition
$\Cal{L} :\, B\to \Cal{G} (H,\alpha ),$ $\Cal{L}(b):=
\Cal{L}^{\prime\prime}(b)\circ\Cal{L}^{\prime}(b).$
A priori this is only a map of sets. The following theorem
shows that we, in fact, get a multiplier.

\smallskip

We will combine the notation of
1.3 and of 2.1, 2.2,
so that, for example, the left
representative of $\Cal{L}^{\prime}(b)$
is $[c_{l,b}^{\prime};\,x^{\prime}_{l,b},h^{\prime}_{l,b},0]$,
the structure form (2.2) of $\Cal{L}^{\prime\prime}$ is
$\langle \,,\rangle^{\prime\prime}$ etc.

\medskip

\proclaim{\quad 3.1.1. Theorem} (a) $\Cal{L}$ is a multiplier.
Its structure form is
$$
\langle b_1,b_2\rangle = \langle b_1,b_2\rangle^{\prime\prime}\,
\langle b_1,b_2\rangle^{\prime}\,
\alpha (h^{\prime}_{l,b_2}, h^{\prime\prime}_{l,b_1})\,
\alpha (h^{\prime}_{l,b_1}, h^{\prime\prime}_{l,b_2}) .
\eqno(3.1)
$$

\smallskip

(b) If $\theta^{\prime}\in \Gamma (\Cal{L}^{\prime})$,
$\theta^{\prime\prime}\in \Gamma (\Cal{L}^{\prime\prime})$,
and $\theta :=\theta^{\prime}\theta^{\prime\prime}$
exists as a formal series, then $\theta\in \Gamma (\Cal{L})$.

\smallskip

(c) If $\alpha$ is unitary and $\Cal{L}^{\prime\prime}$, $\Cal{L}^{\prime}$
are ample, then $\Cal{L}$ is also ample, and all such
pairwise products of theta functions are well defined.  
\endproclaim

\smallskip

{\bf Proof.} (a) From (1.10) we get the left representative of $\Cal{L} (b)$ in the form: 
$$
[c^{\prime\prime}_{l,b}c^{\prime}_{l,b}\,\alpha (h^{\prime}_{l,b},
h^{\prime\prime}_{l,b});\,x^{\prime}_{l,b}A^2_{h_{l,b}^{\prime\prime}},
h^{\prime}_{l,b}+h^{\prime\prime}_{l,b},0].
\eqno(3.2)
$$
We have to check for these representatives the conditions (2.2)
and (2.3) of Lemma 2.2. With (3.2) as input, the right hand side
of (2.2) becomes
$$
h^{\prime}_{l,b_2}(x^{\prime}_{l,b_1} A^2_{h^{\prime\prime}_{l,b_1}})\,
h^{\prime\prime}_{l,b_2}(x^{\prime}_{l,b_1} A^2_{h^{\prime\prime}_{l,b_1}})\,
\times
$$
$$
\alpha (h^{\prime}_{l,b_1},h^{\prime}_{l,b_2})\,
\alpha (h^{\prime\prime}_{l,b_1},h^{\prime\prime}_{l,b_2})\,
\alpha (h^{\prime}_{l,b_1},h^{\prime\prime}_{l,b_2})\,
\alpha (h^{\prime\prime}_{l,b_1},h^{\prime}_{l,b_2}) .
\eqno(3.3)
$$
Using the composability condition
$x^{\prime\prime}_{r,b}=x^{\prime}_{l,b}$ and  (2.1)
in the form  $x^{\prime\prime}_{r,b}=x^{\prime\prime}_{l,b}
A^{-2}_{h_{l,b}^{\prime\prime}}$,
we can replace
the second factor in (3.3) by
$h^{\prime\prime}_{l,b_2}(x^{\prime\prime}_{l,b_1}).$ 
Taking into account (2.2) for the two multipliers,
we obtain the symmetric form (3.1).

\smallskip

It remains to check (2.3). Again from (3.2) we obtain
$$
\frac{c_{l,b_1+b_2}}{c_{l,b_1}c_{l,b_2}}=
\frac{c_{l,b_1+b_2}^{\prime\prime}}{c_{l,b_1}^{\prime\prime}
c_{l,b_2}^{\prime\prime}}\cdot 
\frac{c_{l,b_1+b_2}^{\prime}}{c_{l,b_1}^{\prime}
c_{l,b_2}^{\prime}}\cdot 
\alpha (h^{\prime\prime}_{r,b_1},h^{\prime}_{l,b_2})\,
\alpha (h^{\prime\prime}_{r,b_2},h^{\prime\prime}_{l,b_1}) .
$$
Using (2.2) and replacing the right components $h_r$ by the left
ones, we finally get (3.1).

\smallskip

(b) This immediately follows from (1.11).

\smallskip

(c) If $\alpha$ is unitary, then the logarithm of the modulus
of form (3.1) is negative defined, when this holds
for the factors.

\medskip

{\bf 3.2. Categorical interpretation.} Extending (1.12), we can define
two categories in which composition of multipliers and
multiplication of theta functions become the composition
of morphisms. They have common objects: ``period homomorphisms''
$\xi :\, B\to T(H,1)(K)$ ($B$ is assumed to be fixed;
the most important case is that of a lattice of the same rank as $H$).

\smallskip

In the category $\roman{PIC}\,(H,\alpha )$ morphisms are multipliers:
$$
\roman{Hom}_{{}_\roman{PIC}}\, (\xi ,\eta ):=\{\,\Cal{L}\,|\,
x_r(\Cal{L})=\xi,\, x_l(\Cal{L})=\eta \}
\eqno(3.4)
$$
Here, say, $x_r(\Cal{L})$ is the period map
$b\mapsto x_{r,b}$ associated with $\Cal{L}.$

\smallskip

In the category $\roman{Pic}\,(H,\alpha )$ morphisms are
$$
\roman{Hom}_{{}_\roman{Pic}}\, (\xi ,\eta ):=\{\,\oplus^{\prime}\,\Gamma (\Cal{L})\,|\,
x_r(\Cal{L})=\xi,\, x_l(\Cal{L})=\eta \}
\eqno(3.5)
$$
where the summation $\oplus^{\prime}$ is taken only
over ample multipliers. The resulting morphism spaces are of course graded
by ample multipliers, and the composition is 
compatible with this grading.

\medskip

{\bf 3.3. Multipliers with hidden periods.} We will say
that $\Cal{L}$ and its thetas have {\it hidden periods}, if the
image of $B$ lies in the group of double--sided representatives
$\Cal{G}_d (H,\alpha )$ described in 1.3.3. We continue to
assume that $\alpha^2$ is non--degenerate.

\smallskip

Put $\Cal{L}(b)=[c_b;\,1,h_{l,b}, h_{r,b}]$. From the
condition $\Cal{L}(b_1)\Cal{L}(b_2)=\Cal{L}(b_1+b_2)$ we see that
such a family defines a multiplier iff the alternate form on $B$
$$
\{b_1,b_2\} :=\alpha ( h_{l,b_1},h_{l,b_2})\,
\alpha ( h_{r,b_1},h_{r,b_2})^{-1}
\eqno(3.6)
$$
is symmetric, that is, takes values $\{\pm 1\}$, and moreover
$$
c_{b_1+b_2}=c_{b_1}c_{b_2}\,\{b_1,b_2\}
\eqno(3.7)
$$
for all $b_1,b_2\in B$.

\smallskip

A large class of such multipliers corresponds
to morphisms between quantum tori which are quotients of 
$T(H,\alpha)$, in particular, to the automorphisms
of $T(H,\alpha )$. 

\smallskip

More precisely,
assume that $\Cal{L}$ has the following property:
{\it for any $b$, $h_{r,b}$ is uniquely defined by $h_{l,b}$}.
This holds, for example, if $B\to H:\,b\mapsto h_{l,b}$
has trivial kernel. Put $T_{l,r} :=T(h_{l,r}(B),\alpha )$.
There is a group morphism $f: h_{l}(B)\to h_{r}(B)$
such that $h_{r,b}=f(h_{l,b})$ for all $b$.
If the form (3.6) takes only values $\pm 1$, $f$
is compatible with $\alpha^2$ and therefore 
determines a morphism of quantum tori
$F:\,T_r\to T_l:$ $F^*(e(h))=e(f(h)).$
Conversely, any such morphism of quantum tori
associated with subgroups of $H$ produces
a family of multip[liers with hidden periods
$$
\Cal{L}(b)=[\chi (b);\,1, h_{l,b},f(h_{l,b})]
\eqno(3.8)
$$
where $\chi :\,H\to K^*$ is an arbitrary character.

\smallskip

On an appropriate subset of these data we can also
exhibit a nice functorial description of $\circ$--multiplication:

\smallskip

\proclaim{\quad 3.3.1. Proposition} Consider the following
category $\Cal{T}$: objects are subgroups of maximal rank
in $H$, morphisms are isomorphisms compatible
with $\alpha^2$. 

\smallskip

For each morphism, construct a multiplier as in (3.8).
Then the product of morphisms corresponds to the 
$\circ$--product of multipliers.  
\endproclaim

{\bf 3.3.2. Example.} Consider the two--dimensional torus
$T_q$ described in 0.7: $H=\bold{Z}^2=\bold{Z}h_1\oplus\bold{Z}h_2$,
$\alpha (h_1,h_2)=q,$  $\alpha (h_1,h_1)=
\alpha (h_2,h_2)=1.$ Put $u=e(h_1), v=e(h_2).$
The functional equation for $\theta_q(t)$,
$qu \theta_q (q^2t)=\theta_q (t)$, gives
for two liftings
of $\theta_q$ to $T_q$
$$
quv^{-1}\theta_q(u)v=\theta_q(u),\ u\theta_q(u)u^{-1}=\theta_q(u)
$$
and
$$
qvu\theta_q(v)u^{-1}=\theta_q(v),\ v^{-1}\theta_q(v)v=\theta_q(v).
$$
Hence they are thetas on $T_q$ with hidden periods
corresponding to the the automorphisms
$f:\,h_1\mapsto h_1, h_2\mapsto h_1+h_2,$
$g:\,h_1\mapsto h_1-h_2, h_2\mapsto h_2.$

\smallskip

Cross--multiplying, we get similar functional equations
for $\theta_q(u)\theta_q(v)$:
$$
quv^{-1}\theta_q(u)\theta_q(v)v = \theta_q(u)\theta_q(v),\ 
q^{-1}u\theta_q(u)\theta_q(v)vu^{-1} = \theta_q(u)\theta_q(v).
$$
\medskip

{\bf 3.3.3. Example: Weinstein's theta distribution.} In [We],
A.~Weinstein considers a real quantum torus $T(\bold{Z}^d,\alpha )$
where $\alpha (e(m),e(n))=e^{-\pi i\hbar P(m,n)}$, $P$ being a real
skew--symmetric pairing. He represents the quantum
multiplication by an integral operator with the distribution kernel
$K(x,y,z)=L(y-x,z-x)$ where in his notation 
$$
L(y,z)=\sum_{m,n\in\bold{Z}^d}
e^{-\pi i\hbar P(m,n)-2\pi i(my+nz)} .
\eqno(3.9)
$$
In our general setup, this means encoding
the multiplication on $T(H,\alpha )$ in terms of the
formal theta function on $T(H\oplus H,1)$
$$
\theta_W :=\sum_{(g,h)\in H\times H}\alpha (g,h)\,e(g,h) .
\eqno(3.10)
$$
We can rewrite (3.10) in the most natural way passing
to the torus $T(H\oplus H,\beta )$ where $\beta ((g,0),(0,h)):=
\alpha (g,h)$ and $H\oplus \{0\}, \{0\} \oplus H$
are $\beta$--isotropic. Then we have
$$
u_{1,\beta}(\theta_W)=\sum_{g\in H} e(g,0)\,*_{\beta}\,
\sum_{h\in H}e(0,h)
\eqno(3.11)
$$
where $*_{\beta}$ denotes the multiplication in $T(H\oplus H,\beta )$. 
Each factor in (3.11) is of course the delta distribution at zero
on $T(H,1)$ lifted to $T(H\oplus H,\beta )$. Clearly,
$u_{1,\beta}(\theta_W)$ has the hidden period group
$$
H\oplus H\to\Cal{G} (H\oplus H,\beta ):\ (g,h)\mapsto [1;1,g,-h]_{\beta} .
\eqno(3.12)
$$
Weinstein treats $\theta_W$ as a theta--distribution on the commutative
torus. In order to calculate its automorphy factors,
we can introduce {\it a symmetric pairing} $\langle\,,\rangle$
on $H\oplus H$ such that $\langle (g,h),(g,h)\rangle =\alpha (g,h).$ 
Then $\theta_W=\sum_{k\in H\oplus H} \langle k,k\rangle e(k)$
and moreover, for all $k\in H$
$$
\langle k,k\rangle\, e(k)\,x_k^*(\theta_W)=\theta_W ,
\eqno(3.13)
$$
where $x_k^*(e(j))=\langle k,j\rangle^2.$

\bigskip

\centerline{\bf \S 4. Small Heisenberg groups}

\medskip

{\bf 4.1. Small Heisenberg groups $\Cal{G} (\Cal{L})$.} Consider a
multiplier $\Cal{L}:\, B\to \Cal{G} (H,\alpha )$. Denote
by $\widetilde{\Cal{G}}(\Cal{L})$ the normalizer
of $\Cal{L} (B)$ in $\Cal{G} (H,\alpha )$.

\smallskip

We will call $\Cal{G} (\Cal{L}) := \widetilde{\Cal{G}}(\Cal{L})/\Cal{L}(B)$
{\it the small Heisenberg group} of $\Cal{L}$. Clearly,
it acts upon $\Gamma (\Cal{L})$.
\smallskip

In the following theorem we use the same notation for $\Cal{L}(b)$
as in 2.1. 

\medskip

\proclaim{\quad 4.1.1. Theorem} (a) Assume that
the natural map $\Cal{L} (B)\to T(H,1)(K)\times H$ is
injective. Then
$$
\widetilde{\Cal{G}}_l (\Cal{L}) = \{\,[c;\,\xi ,\gamma ,0 ]\,|\,
\forall b\in B,\ h^-_b(\xi )=\gamma (x_{l,b})\, \alpha^2(h^-_b,\gamma )\,\}.
\eqno(4.1)
$$

\smallskip

(b) Assume that $K$ is large enough in the following sense:
the natural homomorphism $T(H,1)(K)\to T(h^-(B),1)(K)$ is surjective
and all torsion points of its kernel are rational over $K$.
Assume moreover that $\roman{dim}\,\Gamma (\Cal{L})= [H:h^-(B)]
< \infty .$ In this case the representation of $\Cal{G} (\Cal{L})$
in $\Gamma (\Cal{L})$ is irreducible.
\endproclaim

{\bf Proof.} (a) A straightforward calculation gives:
$$
[c;\,\xi ,\gamma ,0]^{-1}[c_{l,b};\,x_{l,b},h^-_b,0]\,
[c;\,\xi ,\gamma ,0]=
$$
$$
[c_{l,b}\,\gamma (x_{l,b})\, \alpha^2(h^-_b,\gamma )h^-_b(\xi );\,x_{l,b},h^-_b
,0] .
\eqno(4.2)
$$
If $[c;\,\xi ,\gamma ,0]$ normalizes $\Cal{L}(B)$
and $\Cal{L} (B)\to T(H,1)(K)\times H$ is
injective, the right hand side of (4.2) must coincide
with $[c_{l,b};\,x_{l,b},h^-_b,0]$. This proves (a).

\smallskip

(b) Let us fix an arbitrary $\gamma\in H$ and consider
the conditions in (4.1) as equations for $\xi$:
$$
\forall b\in B,\ h^-_b(\xi )=\gamma (b)\, \alpha^2(h_b^-
,\gamma )
\eqno(4.3)
$$
These equations uniquely determine the image of this
point in $T(h^-(B),1)(K)$, and we can take for $\xi$
any lift of this image to $T(H,1) (K)$. If the assumptions
of (b) hold, the set of such lifts is non--empty and is
a torsor over the kernel $\Cal{K}$ of    
$T(H,1)(K)\to T(h^-(B),1)(K)$, so that we obtain
an isomorphism $\Cal{G}(\Cal{L})/K^*=\Cal{K}\times H/h^-(B)$
induced by $[c;\,\xi ,\gamma ,0]\mapsto (\xi ,\gamma ).$
The pairing $(\xi ,\gamma )\mapsto \gamma (\xi )$
induces a duality between the finite groups $\Cal{K}$ and $H/h^-(B)$
with values in the group of roots of unity.

\smallskip

Let us now show that $\Cal{G} (\Cal{L})$ and its
representation in $\Gamma (\Cal{L})$  
have exactly the same structure as their
counterparts described in \S 1 of [Mu].

\smallskip

The easiest way to see this is to use the twisting
isomorphism $u_{\alpha ,1}=u_{1,\alpha}^{-1}$ (see (1.6), (1.7))
in order to reduce the situation to the case $\alpha\equiv 1$
which we will assume in the rest of the proof.
Since we assumed that  $\roman{dim}\,\Gamma (\Cal{L})= [H:h^-(B)]
< \infty ,$ the space $\Gamma (\Cal{L})$ splits into
the direct sum of one--dimensional subspaces $\Gamma_{\chi}(\Cal{L})$
labeled by all $\chi \in H/h^-(B)$: each  $\Gamma_{\chi}(\Cal{L})$
consists of thetas whose coefficients $a_h$ vanish outside
the coset class $\chi$ (cf. the proof of Theorem 2.4.1).
Clearly, this splitting coincides with the
splitting into irreducible representations
of $\Cal{K}$ lifted to $\Cal{G}(\Cal{L})$ by $\xi\mapsto\xi^*$.

\smallskip

On the other hand, the set of these subspaces forms
a torsor over any lifting of the group $H/h^-(B)$ to
$\Cal{G}(\Cal{L}).$ More precisely, from (4.1)
with $\alpha \equiv 1$ one sees that the map
$[c;\xi ,\gamma ]\mapsto \gamma$ descends to a well
defined map $\Cal{G}(\Cal{L})\to H/h^-(B)$ which we want
to lift.  

\smallskip

In order to see that such a lifting exists
one can imitate Mumford's reasoning. We omit the details.

\smallskip

The most important consequence of this formalism
is the existence of almost canonical 
bases in all spaces of theta functions $\Gamma (\Cal{L})$
corresponding to (formally) ample multip[liers.
More precisely, the basis depends on the choice
of a lifting of $H/h^-(B)\to \Cal{G}(\Cal{L})$,
it is indexed by the characters of this group, and
defined up to a common constant factor. 

\medskip

{\bf 4.2. Rigidity of comultiplication.} Let $\Cal{A}$
be an abelian variety together with its
parametrization by a torus. Mumford's study
of the morphisms $M^*:\,\Gamma (\Cal{A}\times\Cal{A},
L\boxtimes L)\to \Gamma (\Cal{A}\times\Cal{A},
L^2\boxtimes L^2)$ for symmetric ample $L$ 
admits a straightforward rewriting in terms
of our realization of theta functions.
We can then apply twistings $u_{1,\alpha}$
(resp. $u_{1,\alpha^2}$) to the source (resp. target)
of the morphism $M^*$ (see (1.15)). All formulas not involving
quantum multiplications of thetas will remain the same.

\bigskip

\centerline{\bf References}

\medskip

[AsSch] A.~Astashkevich, A.~S.~Schwarz. {\it Projective modules over
non--commutative tori: classification of modules with
constant curvature connection.} Preprint QA/9904139.

\smallskip

[BEG] V.~Baranovsky, S.~Evens, V.~Ginzburg. {\it Representations
of quantum tori and double--affine Hecke algebras.} Preprint
math.RT/0005024

\smallskip

[Co] A.~Connes. {\it Noncommutative geometry.} Academic Press,
1994.

\smallskip

[CoDSch] A.~Connes, M.~Douglas, A.~Schwarz. {\it Noncommutative
geometry and Matrix theory: compactification on tori.}
Journ. of High Energy Physics, 2 (1998).

\smallskip

[FK] L.~D.~Faddeev, R.~M.~Kashaev. {\it Quantum dilogarithm}.
Preprint hep-th/9310070.

\smallskip

[FKV] L.~D.~Faddeev, R.~M.~Kashaev, A.~Yu.~Volkov. {\it 
Strongly coupled quantum discrete Liouville theory. I: Algebraic
approach and duality.} Preprint hep-th/0006156.

\smallskip

[FV] L.~D.~Faddeev, A.~Yu.~Volkov. {\it Abelian current algebra
and the Virasoro algebra on the lattice.} Preprint
HU--TFT--93--29, 1993.

\smallskip

[KaO] A.~Kapustin, D.~Orlov. {\it Vertex algebras,
mirror symmetry, and D--branes: the case of complex tori.}
preprint hep-th/0010293 .

\smallskip

[KoSo] M.~Kontsevich, Y.~Soibelman. {\it Homological
mirror symmetry and torus fibrations.} Preprint
math.SG/0011041

\smallskip

[Ma1] Yu.~Manin. {\it  Quantized theta--functions.} In: Common
Trends in Mathematics and Quantum Field Theories (Kyoto, 1990), 
Progress of Theor. Phys. Supplement, 102 (1990), 219--228.

\smallskip

[Ma2] Yu.~Manin. {\it Mirror symmetry and quantization of abelian varieties.}
Preprint  math.AG/0005143

\smallskip

[Mu] D.~Mumford. {\it On the equations defining abelian
varieties I.} Inv. Math. 1 (1966), 287 --354.

\smallskip

[Ri] M.~A.~Rieffel. {\it Non--commutative tori --- a case
study of non--commutative differential manifolds.}
In: Cont.~Math., 105 (1990), 191--211.

\smallskip

[RiSch] M.~A.~Rieffel, A.~Schwarz. {\it Morita equivalence
of multidimensional tori.} Preprint math.QA/9803057

\smallskip

[So] Y.~Soibelman. {\it Quantum tori, mirror symmetry
and deformation theory.} Preprint, 2000.

\smallskip

[We] A.~Weinstein. {\it Classical theta functions and quantum tori.}
Publ. RIMS, Kyoto Univ., 30 (1994), 327--333.

\smallskip

[Wi] E.~Witten. {\it Overview of $K$--theory applied to strings.}
Preprint hep-th/0007175

\bigskip

\enddocument